\documentclass[graybox]{svmult}

\usepackage{mathptmx}       
\usepackage{helvet}         
\usepackage{courier}        
\usepackage{type1cm}        
\usepackage{makeidx}         
\usepackage{graphicx}        
\usepackage{multicol}        
\usepackage[bottom]{footmisc}
\usepackage{amsmath,amssymb}

\begin{document}
	
\title*{Bannai--Ito algebras \\and the $osp(1,2)$ superalgebra }
\titlerunning{Bannai--Ito algebras and $osp(1,2)$}
\author{Hendrik De Bie\textsuperscript{1}, Vincent X. Genest\textsuperscript{2}, Wouter van de Vijver\textsuperscript{1}, \and Luc Vinet\textsuperscript{3}}
\authorrunning{De Bie, Genest, Van de Vijver, and Vinet}
\institute{
	\textsuperscript{1}Department of Mathematical Analysis, Faculty of Engineering and Architecture, Ghent University, Ghent, Belgium \email{Hendrik.DeBie@UGhent.be} \email{Wouter.VandeVijver@UGhent.be}
	\and
	\textsuperscript{2}Department of Mathematics, Masschusetts Institute of Technology, Cambridge, USA \email{vxgenest@mit.edu}
	\and
	\textsuperscript{3}Centre de Recherches Math\'ematiques, Universit\'e de Montr\'eal, Montr\'eal, Canada \email{vinet@crm.umontreal.ca}
	}
\maketitle

\abstract{The Bannai--Ito algebra $B(n)$ of rank $(n-2)$ is defined as the algebra generated by the Casimir operators arising in the $n$-fold tensor product of the $osp(1,2)$ superalgebra. The structure relations are presented and representations in bases determined by maximal Abelian subalgebras are discussed. Comments on realizations as symmetry algebras of physical models are offered.}

\section{Introduction}
The Bannai--Ito (BI) algebra $B(3)$ of rank one is the associative algebra with three generators $\Gamma_{12}, \Gamma_{13}, \Gamma_{23}$ obeying the relations
\begin{align}
\label{B3}
	\{\Gamma_{12}, \Gamma_{23}\} = \Gamma_{13} + \omega_{13}, \quad
	\{\Gamma_{12}, \Gamma_{13}\} = \Gamma_{23} + \omega_{23},\quad
	\{\Gamma_{13}, \Gamma_{23}\} = \Gamma_{12} + \omega_{12},
\end{align}	
where $\{A, B\} = AB + BA$ and $\omega_{12}, \omega_{13}, \omega_{23}$ are central. It has been introduced in \cite{Tsujioto&Vinet&Zhedanov-2012} to encode the bispectrality of the BI polynomials. Indeed, the Dunkl shift operators of which the BI polynomials are eigenfunctions, the spectrum variable of the recurrence relation, and the anticommutator of these two operators satisfy \eqref{B3}, up to affine transformations. This algebra is also isomorphic to a degeneration of the double affine Hecke algebra of type $(C_1^{\vee}, C_1)$ \cite{Genest&Vinet&Zhedanov-2016} and has appeared in a variety of contexts. For a review, see \cite{DeBie&Genest&Tsujimoto&Vinet&Zhedanov-2015}.

Interestingly, the algebra \eqref{B3} arises in the context of the representation theory of the Lie superalgebra $osp(1,2)$, more specifically in the recoupling schemes for the tensor product of three irreducible representations. In this framework, the BI polynomials are seen to be essentially the Racah coefficients of $osp(1,2)$, that is the elements of the matrices relating the bases associated to the coupling of the first two factors of the three-fold product to the basis corresponding to the situation where the last two factors are initially regrouped. This connection of $B(3)$ to $osp(1,2)$ extends to $n$-fold tensor products and leads to the BI algebra $B(n)$ of arbitrary rank.

This will be presented in the following. We shall also give indications of how representations of $B(n)$ can be constructed in bases associated to maximal Abelian subalgebras. We shall conclude by mentioning some applications of $B(n)$.

\section{$osp(1,2)$ and the Bannai--Ito algebra}

The $osp(1,2)$ superalgebra can be presented as follows. It is generated by two odd elements $J_{\pm}$ and one even element $J_0$ that obey
\begin{align*}
[J_0, J_{\pm}] = \pm J_{\pm}, \quad \{J_{+}, J_{-}\} = 2J_0,
\end{align*}
with $[a, b] = ab - ba$. The $\mathbb{Z}_2$ grading can be accounted for by introducing the grade involution $P$ and including the relations
\begin{align*}
[J_0, P] = 0,\qquad \{J_{\pm}, P\} = 0,\qquad P^2 = 1.
\end{align*}
The Casimir operator in the universal enveloping algebra $\mathcal{U}(osp(1,2))$ 
\begin{align}
\label{Casimir}
\Gamma = \frac{1}{2} \left([J_{-}, J_{+}] - 1\right)P = J_0 P -J_{+}J_{-}P - P/2,
\end{align}
is found to commute with all generators. There is an algebra morphism $\Delta: osp(1,2)\rightarrow osp(1,2)\otimes osp(1,2)$ called comultiplication that acts as follows on the generators:
\begin{align*}
\Delta(J_0) = J_0 \otimes 1 + 1 \otimes J_0,\qquad \Delta(J_{\pm}) = J_{\pm} \otimes P + 1 \otimes J_{\pm},\qquad \Delta(P) = P \otimes P,
\end{align*}
and is coassociative, i.e. $(\Delta \otimes 1)\Delta = \Delta(1\otimes \Delta)$. The coproduct can be iterated to form higher tensor powers of $osp(1,2)$. For a positive integer $n$, define $\Delta^{(n)}: osp(1,2)\rightarrow osp(1,2)^{\otimes n}$ as $\Delta^{(n)} = (1^{\otimes (n-2)}\otimes \Delta)\circ \Delta^{(d-1)}$, with $\Delta^{(1)} = \mathrm{Id}$. 

Let $[n] = \{1,\ldots, n\}$ and let $A = \{a_1,\ldots, a_{k}\}$ be an ordered $k$-subset of $[n]$. For $1\leq k \leq n$, one has a realization of $osp(1,2)$ in $osp(1,2)^{\otimes n}$ for any $A$. This realization, denoted by $osp^{A}(1,2)$, has generators
\begin{align*}
J_{\pm}^{A} = \sum_{a_i\in A} J_{\pm}^{(a_i)}\prod_{j=a_i+1}^{a_k}P^{(j)},\qquad J_0^{S}= \sum_{a_i\in A} J_0^{(a_i)}, \qquad P^{A} = \prod_{a_i\in A} P^{(a_i)},
\end{align*} 
where $J_{\pm}^{(i)}$, $J_0^{(i)}$, $P^{(i)}$ denote the generators of the $i$\textsuperscript{th} factor of $osp(1,2)$ in $osp(1,2)^{\otimes n}$. We can now define the following elements in $\mathcal{U}(osp(1,2)^{\otimes n})$:
\begin{align}
\label{Casimir-A}
\Gamma_{A} =J_0^{A} P^{A} -J_{+}^{A}J_{-}^{A}P^{A} - P^{A}/2.
\end{align}
Clearly $\Gamma_{\{i\}}$, $i=1,\ldots, n$, are the Casimir elements corresponding to each of the factors in $osp(1,2)^{\otimes n}$; these are constant multiples, say $\lambda_i$, if one considers products of irreducible representations. $\Gamma^{[n]}$ is the total Casimir operator of $osp(1,2)^{\otimes n}$. It will be convenient to take $\Gamma_{\emptyset} = -1/2$. We define the Bannai--Ito algebra $B(n)$ as the algebra generated by the elements $\Gamma_{A}$ with $A\subset [n]$.

Let us now determine the structure relations. Consider first the case $n = 3$. There are seven generators in this instance: $\Gamma_{\{i\}} \equiv \Gamma_i$, $i = 1,2,3$, $\Gamma_{\{i,j\}} \equiv \Gamma_{ij}$ for $\{i,j\} = \{1,2\}, \{1,3\}, \{2,3\}$, and $\Gamma_{[3]} \equiv \Gamma_{123}$. Here, $\Gamma_1$, $\Gamma_2$, $\Gamma_3$ and $\Gamma_{123}$ are central. A direct calculation shows that
\begin{align*}
\{\Gamma_{ij}, \Gamma_{jk}\} = \Gamma_{ik} + 2 \Gamma_{j} \Gamma_{ijk} + 2 \Gamma_{i} \Gamma_{k}, \qquad i \neq j \neq k,
\end{align*}
which coincides with the $B(3)$ defining relations \eqref{B3} when the central $\omega_{ik}$ are identified as $\omega_{ik} = 2 \Gamma_{j}\Gamma_{ijk} + 2\Gamma_{i}\Gamma_{k}$. This is the result obtained in \cite{Genest&Vinet&Zhedanov-2014}. The structure relations for the higher rank extension are obtained from this result through the following argument. Take any triple of pairwise disjoint subsets of $[n]$ called $K$, $L$, and $M$. There is an obvious isomorphism
\begin{align*}
osp^{K}(1,2) \otimes osp^{L}(1,2) \otimes osp^{M}(1,2) \cong osp(1,2) \otimes osp(1,2) \otimes osp(1,2),
\end{align*}
which leads to an embedding of $B(3)$ into $B(n)$. Indeed, in view of this isomorphism, the Casimir elements $\Gamma_{K}$, $\Gamma_{L}$, $\Gamma_{M}$, $\Gamma_{K \cup L}$, $\Gamma_{K \cup M}$, $\Gamma_{L \cup M}$, and $\Gamma_{K \cup L \cup M}$ will generate $B(3)$ and we shall have for instance 
\begin{align}
\label{Above}
\{\Gamma_{K \cup L}, \Gamma_{L \cup M}\} = \Gamma_{K \cup M} + 2 \Gamma_{L} \Gamma_{K \cup L \cup M} + 2 \Gamma_{K} \Gamma_{M}.
\end{align}
We wish to know $\{\Gamma_{A}, \Gamma_{B}\}$ for any two subsets $A$ and $B$ of $[n]$. To that end, take $K = A\setminus B$, $L = A \cap B$, $M = B\setminus A$ and make the corresponding Casimir operators in $\mathcal{U}(osp(1,2)^{\otimes n})$ using \eqref{Casimir-A}. The relation \eqref{Above} becomes
\begin{align}
\label{B(n)}
\{\Gamma_{A}, \Gamma_{B}\} = \Gamma_{(A\cup B)\setminus (A\cap B)} + 2\Gamma_{A\cap B}\Gamma_{A\cup B} + 2 \Gamma_{A\setminus (A\cap B)} \Gamma_{B\setminus (A\cap B)}.
\end{align}
This provides the desired structure relations for $B(n)$, namely the relations obeyed by the Casimir elements $\Gamma_{A}$ labeled by subsets $A$ of $[n]$ \cite{DeBie&Genest&Vinet-2016}.

\section{Maximal Abelian subalgebras, representation bases, and connection coefficients}
We wish to indicate here how representations of $B(n)$ can be obtained from the knowledge of representations of $B(3)$. To that end, we shall first introduce bases for representation spaces that are associated to maximal Abelian subalgebras of $B(n)$.

\subsection{Maximal Abelian subalgebras}
We readily see from \eqref{B(n)} that $[\Gamma_A, \Gamma_B] = 0$ if $A\subset B$, $B\subset A$ or $A\cap B = \emptyset$; recall that $\Gamma_{\emptyset} = -1/2$. It follows that $\mathcal{Y}_{n} = \langle \Gamma_{[2]}, \Gamma_{[3]}, \ldots, \Gamma_{[n-1]}\rangle$ forms an Abelian subalgebra (AS) of $B(n)$ that is readily seen to be maximal. Note that $\Gamma_{\emptyset}$, $\Gamma_{[1]}$ and $\Gamma_{[n]}$ are not included in $\mathcal{Y}_n$ as they are central in $B(n)$. Other such maximal AS can be obtained by applying a permutation and taking $\pi \mathcal{Y}_n = \langle \Gamma_{\pi[2]}, \Gamma_{\pi [3]},\ldots, \Gamma_{\pi [n-1]}\rangle$ 

\subsection{Bases for representation spaces}
Bases for representations spaces can now be obtained by taking their elements to be joint eigenvectors of the operators (hereafter denoted by the same symbols) representing the generators of the various maximal AS. Given one such basis, one would wish to provide the action of the generators in the complement of the AS in order to construct the representation of $B(n)$. We shall indicate how this can be accomplished from the knowledge of the connection coefficients between bases associated to different maximal AS. With this understood, we shall complete the picture with a characterization of the connection coefficients. 

Suppose that a basis has been picked and that we want to give the action of a generator $\Gamma$ on the elements of this basis. It is easy to see that every generator of $B(n)$ belongs to a maximal AS. There is a thus another basis, call it prime, in which $\Gamma$ is diagonal. Now if the connection coefficients between the elements of the original bases and those of the prime basis are known, it follows from linear algebra that the action of $\Gamma$ in the original basis can be written down. This applies to any generator. Hence if all bases associated to maximal AS can be connected, the action of all generators in a single basis can be obtained with the help of the connection coefficients.

\subsection{Connection coefficients}
Irreducible representations of $B(3)$ have been constructed and as a result the connection coefficients (CCs) between the bases associated to the AS generated respectively by $\Gamma_{12}$, $\Gamma_{13}$ and $\Gamma_{23}$ are known; see \cite{DeBie&Genest&Vinet-2016, DeBie&Genest&Vinet-2016-2,Genest&Vinet&Zhedanov-2014, Tsujioto&Vinet&Zhedanov-2012}. We shall simply set the notation and recall the main features. Consider an irreducible representation of $B(3)$ and let $\langle \phi_k\rangle$ be a set of basis vectors on which $\Gamma_{12}$ acts diagonally, say $\Gamma_{12}\phi_k = \mu_k \phi_k$. The central elements are multiples of the identity: $\Gamma_i \phi_k = \lambda_i \phi_k$, $\Gamma_{123} \phi_k = \lambda_{123} \phi_k$ for $i=1,2,3$. It is found that $\Gamma_{13}$ and $\Gamma_{23}$ act in a tridiagonal fashion in the basis $\langle \phi_k \rangle$ and one has for instance $\Gamma_{23} \phi_k = a_{k,k-1} \phi_{k-1} + a_{k,k} \phi_{k} + a_{k,k+1}\phi_{k+1}$. The coefficients $a_{k,k}$, $a_{k,k\pm 1}$ have been explicitly determined from the properties of $B(3)$. Now if $\langle \psi_k\rangle $ denotes the basis in which $\Gamma_{23}$ is diagonal, the CCs defined by
\begin{align*}
\psi_k = \sum_{s} B_{ks}(\lambda_1, \lambda_2, \lambda_3, \lambda_{123}) \phi_{s},
\end{align*}
are nothing else than the Racah coefficients of $osp(1,2)$. Knowing the action of $\Gamma_{23}$ in both bases $\langle \phi_k \rangle$ and $\langle \psi_k \rangle$, one finds that $B_{ks}$ satisfies a three-term recurrence relation which shows that these CCs can be expressed in terms of BI polynomials.

Let us now  discuss the rank two case $B(4)$ to illustrate how one bootstraps from the rank one to higher ranks. First consider the CCs between two bases associated to two maximal AS that differ by only one generator. An example is $(\Gamma_{12}, \Gamma_{123})$, $(\Gamma_{12},\Gamma_{124})$. $\Gamma_{123}$ and $\Gamma_{124}$ will preserve the eigenspaces of the common generator $\Gamma_{12}$. Now note that $\Gamma_{123}$ and $\Gamma_{124}$ are also generators of a rank one BI algebra. Indeed, let $K = \{1,2\}$, $L = \{3\}$, $M = \{4\}$, $\Gamma_{K \cup L} = \Gamma_{123}$, $\Gamma_{K \cup M} = \Gamma_{124}$, $\Gamma_{L \cup M} = \Gamma_{34}$ provide an embedding of $B(3)$ into $B(4)$. These generators all commute with $\Gamma_{12}$ and the basis vectors with fixed eigenvalues of $\Gamma_{12}$ will support representations of $B(3)$. The representation theory of the rank one BI algebra tells us that the CCs will again be BI polynomials. Thus is $\langle \phi_{j_1, j_2}\rangle$ and $\langle \psi_{j_1,j_2}\rangle$ are the bases diagonalizing the maximal AS of our example with $\Gamma_{12} \phi_{j_1, j_2} = \mu_{j_1}^{12}\phi_{j_1,j_2}$ and $\Gamma_{12}\psi_{j_1,j_2}=\mu_{j_1,j_2}^{12}\psi_{j_1,j_2}$ we have
\begin{align*}
\psi_{j_1,j_2} = \sum_{k} W_{j_2 k}(\mu_{j_1}^{12}, \lambda_3, \lambda_4, \lambda_{1234}) \phi_{j_1,k}.
\end{align*}
This then allows to obtain the actions of $\Gamma_{124}$ and $\Gamma_{34}$ on the basis vectors $\phi_{j_1,j_2}$. To find the action of other generators, one must consider the relations of the basis $\langle \phi_{j_1,j_2}\rangle$ with other subalgebra-type bases. It can be seen that there is always a path between any given basis to all the others that is made out of intermediary segments where the corresponding AS only differ by one element. The CCs between any two bases are then obtained by iterating for each of those segments the procedure just described when only one generator is different. The resulting CCs will hence be given by a product of BI polynomials. Furthermore, since all generators are part of maximal AS and are diagonal in the corresponding bases, the knowledge of the CCs allows to obtain the action of all generators in a chosen basis. These considerations extend from $B(4)$ to $B(n)$ and it follows that the representations of $B(n)$ in a fixed subalgebra-type basis can be fully characterized.

\section{Conclusion}
We conclude by mentioning that the Bannai--Ito algebra $B(n)$ has arisen in various systems. These models are obtained from particular realizations of $osp(1,2)$:
\begin{itemize}
	\item The Dunkl Laplacian when Dunkl operators are used \cite{Genest&Vinet&Zhedanov-2015-2};
	\item The Dirac--Dunkl equation when Clifford algebras are introduced \cite{DeBie&Genest&Vinet-2016, DeBie&Genest&Vinet-2016-2}
	\item The superintegrable model with reflections with Hamiltonian
	\begin{align*}
	H = \sum_{1 \leq i < j \leq n} J_{ij}^2 + \sum_{k=1}^{n}\frac{\mu_k(\mu_k-r_k)}{x_k^2}
	\end{align*}
	with $r_k f(x_1,\ldots, x_k,\ldots, x_n) = f(x_1,\ldots,-x_k,\ldots, x_n)$ and $J_{ij}= i(x_j \partial_{x_i} - x_{i}\partial_{x_j})$ when gauge-transformed parabosonic operators are called upon \cite{DeBie&Genest&Lemay&Vinet-2016, Genest&Vinet&Zhedanov-2014-2}. In this connection, see \cite{Plyushchay1996, Wigner1950, Yang1951}
\end{itemize}
All these models have the Bannai--Ito algebra as symmetry algebra. We may suspect that $B(n)$ and its representations will keep appearing in different guises.
\section*{Acknowledgments}
This paper was completed during a stay of LV at the School of Mathematical Sciences of the Shanghai Jia Tong University as Chair Visiting Professor. HDB is supported by the Fund for Scientific Research-Flanders (FWO-V), project ``Construction of algebra realizations using Dirac-operators'', grant G.0116.13N. VXG holds a postdoctoral fellowship from the Natural Science and Engineering Research Council of Canada (NSERC). The research of LV is supported in part by NSERC.

\end{document}